\documentclass[a4paper,12pt]{article}
\usepackage{amsmath,amsthm,verbatim,amssymb,amsfonts,amscd, graphicx}
\usepackage[all]{xy}
\usepackage{xcolor}
\usepackage{geometry}
\usepackage{graphicx}
\usepackage{hyperref}
\usepackage{enumitem}
\definecolor{darkblue}{rgb}{0,0,0.3}
\definecolor{urlblue}{rgb}{0,0,0.7}
\hypersetup{
	colorlinks=true,
	citecolor=blue,
	linkcolor=darkblue,
	urlcolor=urlblue,
}

\newcommand{\RR}{\mathbb{R}}

\DeclareMathOperator{\tr}{tr}

\DeclareMathOperator{\diam}{diam}

\DeclareMathOperator{\Hess}{Hess}
\DeclareMathOperator{\Ric}{Ric}

\newcommand{\D}{\nabla}
\newcommand{\p}{\partial}

\newcommand{\metric}[2]{\langle#1\,,\,#2\rangle}

\renewcommand{\bar}{\overline}
\renewcommand{\tilde}{\widetilde}

\newcommand{\updatetag}[2]{}
\newcommand{\forcedlinebreak}[1]{\vspace{#1}}
\newcommand\blfootnote[1]{%
	\begingroup
	\renewcommand\thefootnote{}\footnote{#1}%
	\addtocounter{footnote}{-1}%
	\endgroup
}

\newtheorem{theorem}{Theorem}[section]
\newtheorem{lemma}[theorem]{Lemma}
\newtheorem{cor}[theorem]{Corollary}
\newtheorem{defn}[theorem]{Definition}
\numberwithin{equation}{section}
\theoremstyle{definition}
\newtheorem{remark}[theorem]{Remark}

\newcommand{\C}{\mathcal{C}}
\newcommand{\V}{\mathcal{V}}
\renewcommand{\H}{\mathcal{H}}

\begin{document}
	
\title{Dimension constraints in some problems involving intermediate curvature}
\author{Kai Xu}
\date{}
\maketitle

\begin{abstract}
	In \cite{Brendle-Hirsch-Johne_2022} Brendle-Hirsch-Johne proved that $T^m\times S^{n-m}$ does not admit metrics with positive $m$-intermediate curvature when $n\leq 7$. Chu-Kwong-Lee showed in \cite{Chu-Kwong-Lee_2022} a corresponding rigidity statement when $n\leq 5$. In this paper, we show the sharpness of the dimension constraints by giving concrete counterexamples in $n\geq 7$ and extending the rigidity result to $n=6$. Concerning uniformly positive intermediate curvature, we show that simply-connected manifolds with dimension $\leq 5$ and bi-Ricci curvature $\geq 1$ have finite Urysohn 1-width. Counterexamples are constructed in dimension $\geq 6$.
\end{abstract}

\section{Introduction}

\blfootnote{\textit{Mathematics Subject Classification (2020)}. 53C21; 53C42.}In this paper, we study a notion of curvature that interpolates between Ricci and scalar curvatures. For $\{e_i\}_{i=1}^n$ an orthonormal frame at a point on a manifold $M$, we organize the sectional curvatures $\sec(e_i,e_j)$ in an $n\times n$ matrix. Thus Ricci curvature is the sum of entries in the first row, and scalar curvature is twice the sum of entries in the upper triangular part. The curvature quantity that we study, which we call \textit{$m$-intermediate curvature} throughout this paper following Brendle-Hirsch-Johne \cite{Brendle-Hirsch-Johne_2022}, is defined as the sum of the first $m$ rows of the upper triangular part.

\begin{defn}[see also {\cite[Definition 1.1]{Brendle-Hirsch-Johne_2022}}]\label{def-intro:Cm}
	Let $M$ be an $n$-manifold, and $1\leq m\leq n-1$. Suppose $\{e_1,e_2,\cdots,e_n\}$ is an orthonormal frame at $x\in M$. The $m$-intermediate curvature of $M$, denoted by $\C_m$, is defined by
	\[\C_m(e_1,\cdots,e_m):=\sum_{p=1}^m\sum_{q=p+1}^n R(e_p,e_q,e_q,e_p).\]
	We say that $M$ has $m$-intermediate curvature lower bound $\lambda$ (abbreviated as $\C_m\geq\lambda$), if
	\[\C_m(e_1,\cdots,e_m)\geq\lambda\text{\ \ for any orthonormal frame $\{e_i\}$ at any $x\in M$}.\]
	We say that $M$ has positive $m$-intermediate curvature (\,$\C_m>0$), if $\C_m(e_1,\cdots,e_m)>0$ for any orthonormal frame $\{e_i\}$. When $m=2$ (resp. $m=3$), the $m$-intermediate curvature is also called bi-Ricci (resp. tri-Ricci) curvature.
\end{defn}

Apart from Definition \ref{def-intro:Cm}, there are several other notions of interpolating curvatures in the literature. For the partial sum of eigenvalues of Ricci curvature and curvature operator, see Hirsch-Kazaras-Khuri-Zhang \cite{Hirsch-Kazaras-Khuri-Zhang_2022}, Wolfson \cite{Wolfson_2012} and Bohm-Wilking \cite{Bohm-Wilking_2008}, Petersen-Wink \cite{Petersen-Wink_2021}. For curvature quantities that interpolate between sectional and Ricci curvature, see Shen \cite{Shen_1989}, Wu \cite{Wu_1987}. The notion of intermediate curvature in Definition \ref{def-intro:Cm} suits particularly with the technique of minimal weighted slicing technique developed in \cite{Brendle-Hirsch-Johne_2022}, which enables one to obtain topological obstruction theorems (see Theorem \ref{thm-intro:Brendle-Hirsch-Johne}). We remark that, on the other hand, it remains unknown whether $S^2\times T^2$ admits metrics whose sum of the first two smallest eigenvalues of Ricci curvature is positive.

The following property is useful: $\C_m(e_1,\cdots,e_m)$ only depends on the subspace which $\{e_1,\cdots,e_m\}$ spans, thus is a function on the $m$-Grassmannian bundle over $M$. To see this, we orthogonally split $T_x M=V\oplus V^\perp$ where $V=\text{span}\{e_1,\cdots,e_m\}$. The metric tensor also splits as $g=g_1\oplus g_2$. Then we have
\[\C_m(e_1,\cdots,e_m)=\big\langle R_{ijkl}\,,\,\frac12(g_1)_{il}(g_1)_{jk}+(g_1)_{il}(g_2)_{jk}\big\rangle,\]
proving our claim. If a metric $g$ has $\C_m>0$ at a point, then it has positive scalar curvature at the same point. More properties of intermediate curvature can be found in \cite{Brendle-Hirsch-Johne_2022}.

\vspace{12pt}

An important model space for studying intermediate curvature is the product of flat torus and round sphere $M=T^m\times S^{n-m}$: $M$ has positive $\C_{m+1}$ but only non-negative $\C_m$. It is natural to ask a topological obstruction problem: whether $T^m\times S^{n-m}$ admit metrics with positive $\C_m$. This question extends two well-known results. For the case $m=1$ (i.e. $\Ric>0$), Bonnet-Myers' theorem or Bochner's formula gives negative answer. The case $m=n-1$ (i.e. $M=T^n$ with positive scalar curvature) is Geroch's conjecture, for which the answer is also negative as shown by Schoen-Yau \cite{Schoen-Yau_1979a} \cite{Schoen-Yau_2017} and Gromov-Lawson  \cite{Gromov-Lawson_1983}. The extension of topological obstruction to intermediate curvature is recently proved by Brendle-Hirsch-Johne \cite{Brendle-Hirsch-Johne_2022} up to dimension 7.

\begin{theorem}[\cite{Brendle-Hirsch-Johne_2022}]\label{thm-intro:Brendle-Hirsch-Johne}
	For $n\leq 7$, there is no metric on $T^m\times M^{n-m}$ with $\C_m>0$, where $M^{n-m}$ is any $(n-m)$-dimensional closed manifold.
\end{theorem}

Chu-Kwong-Lee \cite{Chu-Kwong-Lee_2022} later proved a corresponding rigidity statement.

\begin{theorem}[\cite{Chu-Kwong-Lee_2022}]\label{thm-intro:Chu-Kwong-Lee}
	Let $n\leq5$. If a metric $g$ on $T^m\times M^{n-m}$ satisfies $\C_m\geq0$, then $g$ splits isometrically as the product of a flat torus $T^m$ with a metric on $M^{n-m}$ with nonnegative Ricci curvature.
\end{theorem}

In this paper we are focused on the dimension conditions that appear in the above theorems. Our result unexpectedly shows that concrete counterexamples exist in higher dimensions, which are not due to the smoothness issue of area-minimizing hypersurfaces. More precisely we have:

\begin{theorem}\label{thm-intro:countereg_Ricm}
	Assume $m\geq 2$, $n\geq m+2$. Let $F=F(x_1,\cdots,x_m)$ be a positive function on the torus $T^m$. Consider the following metric on $S^{n-m}\times T^m$:
	\begin{equation}\label{eq-intro:countereg_Ricm}
		g=\epsilon^2F^2h+F^{-2\frac{n-m}{m-1}}\sum_{i=1}^m dx_i^2,
	\end{equation}
	where $h$ is the standard metric on $S^{n-m}$. Then we have
	\[\C_m\Big(\frac{\p x_1}{|\p x_1|},\cdots,\frac{\p x_m}{|\p x_m|}\Big) =
	(n-m)\Big[\frac{(n-m)(m-2)}{2(m-1)}-1\Big] F^{2\frac{n-m}{m-1}-2}
		\sum_{i=1}^m\big(\frac{\p F}{\p x_i}\big)^2.\]
	If $n(m-2)\geq m^2-2$ and $\epsilon$ is sufficiently small (depending on the $C^2$ norm of $\log F$), then $g$ satisfies $\C_m\geq0$ everywhere. If further $n(m-2)>m^2-2$, then $\C_m>0$ at all non-critical points of $F$.
\end{theorem}

The algebraic relation $n(m-2)\geq m^2-2$ obtained here is consistent with the algebraic inequalities in \cite[Lemma 3.14]{Brendle-Hirsch-Johne_2022} and \cite[Lemma 2.3]{Chu-Kwong-Lee_2022}. From Theorem \ref{thm-intro:countereg_Ricm}, we obtain metrics on $T^4\times S^4$ and $T^3\times S^5$ with positive intermediate curvature except on the critical points of $F$. The latter issue can be resolved by a small perturbation. As a result, we obtain examples that show the sharpness of the bound $n\leq7$ in Theorem \ref{thm-intro:Brendle-Hirsch-Johne}.

\begin{cor}\label{cor-intro:perturbation}
	Assume and all the conditions in Theorem \ref{thm-intro:countereg_Ricm}, and $n(m-2)>m^2-2$. If $F$ is a Morse function on $T^m$ and $g$ is as in \eqref{eq-intro:countereg_Ricm}, then there exists a small perturbation of $g$ that satisfies $\C_m>0$ strictly.
\end{cor}

Regarding Theorem \ref{thm-intro:Chu-Kwong-Lee}, we obtain non-trivial metrics on $T^3\times S^4$ and $T^4\times S^3$ with non-negative intermediate curvatures, giving counterexamples for $n=7$. The case $n=6$ is not included in Theorem \ref{thm-intro:Chu-Kwong-Lee}. However, the proof in \cite{Chu-Kwong-Lee_2022} can be improved and thus extend to $n=6$. We refer to Section \ref{sec:dim6} for details.

\begin{theorem}\label{thm-intro:rigidity_dim6}
	The statement of Theorem \ref{thm-intro:Chu-Kwong-Lee} holds for $n=6$ as well.
\end{theorem}

Given the formula (\ref{eq-intro:countereg_Ricm}), it is elementary to verify all the conclusions. The important aspect of this result is to actually find (\ref{eq-intro:countereg_Ricm}). To better explain the idea, let us briefly summarize the proof of Theorem \ref{thm-intro:Brendle-Hirsch-Johne} in \cite{Brendle-Hirsch-Johne_2022}. If there exists a metric as stated in the theorem, we find a chain of hypersurfaces $M\supset\Sigma_1\supset\cdots\supset\Sigma_m$, with $\Sigma_1$ being a stable minimal hypersurface in $M$, and each $\Sigma_i$ being a minimizer of a certain weighted area functional on $\Sigma_{i-1}$. Combining all the stability inequalities of $\Sigma_i$, one obtains an inequality of the form
\begin{equation}\label{eq-intro:aux1}
	0\leq\int_{\Sigma_m}\big(-\C_m-\V_1-\V_2-\cdots\big),
\end{equation}
where $\V_1, \V_2,\cdots$ are algebraic expressions involving the second fundamental forms of $\Sigma_i$. We would obtain a contradiction if $\V_i\geq0$ for all $i$. As shown in \cite{Brendle-Hirsch-Johne_2022}, this is true when the inequality $n(m-2)\geq m^2-2$ is satisfied.

Our further observation is that: the validity of the algebraic inequalities $\V_i\geq0$ suggests the existence of counterexamples. If the algebraic inequality $\V_i\geq0$ fails to be always true in a certain dimension, then we reverse engineer a metric that captures the negative part. Apart from this, we manage to achieve equality for all the other inequalities involved. As a result, (\ref{eq-intro:aux1}) will force $\C_m\equiv-(\V_1+\V_2+\cdots)>0$ to hold, and thus give rise to nontrivial counterexamples. This process finally yields (\ref{eq-intro:countereg_Ricm}). A more detailed explanation of this process is contained in subsection \ref{subsec:engi}. \\

The next problem concerns diameter bounds under the presence of uniformly positive intermediate curvatures. By Bonnet-Myers' theorem, manifolds with uniformly positive Ricci curvature have bounded diameter. When trying to generalize this fact to bi-Ricci curvature, one needs to pass to a stable minimal surface.

\begin{theorem}[Shen-Ye \cite{Shen-Ye_1996}]\label{thm-intro:Shen-Ye}
	Let $M$ be a (complete) manifold with dimension $n\leq 5$ and bi-Ricci curvature $\C_2\geq\lambda>0$. Then any two-sided stable minimal hypersurface in $M$ has diameter $\leq C(n)\lambda^{-1/2}$.
\end{theorem}

We first note that, with the stability inequality for minimal hypersurfaces, Theorem \ref{thm-intro:Shen-Ye} can be interpreted into a statement in terms of the stable minimal surface (without referring to the ambient manifold $M$).

\begin{defn}\label{def-intro:min_Ric}
	On a Riemannian manifold $M$, we define the minimum Ricci scalar
	\[R_0(x):=\min_{e\in T_pM,|e|=1}\Ric_M(e,e),\quad x\in M.\]
	Equivalently, $R_0$ is the minimal eigenvalue of Ricci curvature at each point.
\end{defn}

\begin{lemma}\label{lemma-intro:weakly_Ric>0}
	Suppose $M$ satisfies bi-Ricci curvature lower bound $\C_2\geq\lambda$, and $\Sigma\subset M$ is a two-sided stable minimal hypersurface. Then
	\begin{equation}\label{eq-intro:weak_Ricci_lower_bound}
		\lambda_1(-\Delta_\Sigma+R_0)\geq\lambda
	\end{equation}
	(where $\lambda_1$ denotes the first eigenvalue), or equivalently, we have
	\begin{equation}
		\int_\Sigma|\D\varphi|^2+R_0\varphi^2\geq\lambda\int_\Sigma\varphi^2
	\end{equation}
	for all $\varphi\in C^\infty(\Sigma)$.
\end{lemma}

Theorem \ref{thm-intro:Shen-Ye} can now be understood as a Bonnet-Myers' theorem for the spectral Ricci curvature lower bound (\ref{eq-intro:weak_Ricci_lower_bound}). We are concerned with the range of dimensions for which such a theorem holds. Generalizing (\ref{eq-intro:weak_Ricci_lower_bound}) by adding a coefficient in front of $R_0$ helps us see clearly the optimal form of a dimension constraint. We have:

\begin{theorem}\label{thm-intro:weak-Bonnet-Myers}
	Let $\beta>0$, $n\leq 7$. Suppose an $n$-dimensional complete manifold $\Sigma$ satisfies
	\[\lambda_1(-\Delta_\Sigma+\beta R_0)\geq\lambda>0,\]
	then we have $\diam(\Sigma)\leq C(n,\beta)\lambda^{-1/2}$ as long as
	\begin{equation}\label{eq-intro:dim-beta_ineq}
		\left\{\begin{aligned}
			& \beta\geq\frac12\quad(n=3), \\
			& \beta>\frac14(n-1)\quad(n\ne 3).
		\end{aligned}\right.
	\end{equation}
	In particular, $\Sigma$ must be compact. There exist non-compact counterexamples when (\ref{eq-intro:dim-beta_ineq}) is not satisfied.
\end{theorem}

The restriction $n\leq 7$ is due to smoothness of area-minimizing currents. Theorem \ref{thm-intro:Shen-Ye} corresponds to the special case $\beta=1$, which confirms the sharpness of the condition $n\leq 5$ there. A proof of Theorem \ref{thm-intro:weak-Bonnet-Myers} (without counterexamples) was given in Shen-Ye \cite{Shen-Ye_1997} using the stability inequality for weighted minimal geodesics. Moreover, Carron-Rose \cite{Carron-Rose_2021} proved similar diameter bound under a stronger hypothesis $\beta>\frac{n+4}3$, $n\geq 3$. Here we give a new proof of Theorem \ref{thm-intro:weak-Bonnet-Myers} based on the $\mu$-bubble technique. An introduction to $\mu$-bubbles, including references, is included in subsection \ref{subsec:mu-bubbles}. This new proof has the benefit that counterexamples can be reverse-engineered, especially for the critical case $\beta=\frac{n-1}4$ (note that hyperbolic spaces are counterexamples for $\beta<\frac{n-1}4$). The reverse engineering process also gives rise to interesting counterexamples with finite volume for $\beta<\frac{n-1}4$, $n\geq3$, which behaves differently from the hyperbolic space at infinity.

We observe that, the usage of $\mu$-bubbles imposes a stronger dimensional constraint than minimal surfaces. For example, using weighted minimal surfaces one can show that $S^{n-1}\times S^1$ ($n\leq 7$) does not admit metrics with $\lambda_1(-\Delta+\beta R_0)>0$ under a weaker bound $\beta\geq\frac{n-2}{n-1}$ (see Remark \ref{rmk-Ric2>0:top_obstruction}). This threshold has also appeared in Bour-Carron \cite{Bour-Carron_2017}. \\

Finally, we relate intermediate curvature to the notion of macroscopic dimension. A manifold $M$ is said to have \textit{macroscopic dimension} $\leq k$, or to have finite \textit{Urysohn $k$-width}, if there exists a simplicial complex $X$ with dimension $\leq k$, and a continuous map $f:M\to X$, such that each fiber $f^{-1}(p)$ has uniformly bounded diameter. By Bonnet-Myers' theorem, manifolds with $\Ric\geq\lambda>0$ have finite Urysohn 0-width. It is conjectured by Gromov that an $n$-manifold with uniformly positive scalar curvature has finite Urysohn $(n-2)$-width. The dimension 3 case has been answered affirmatively by Gromov-Lawson, Katz, Liokumovich-Maximo \cite{Gromov-Lawson_1983, Katz_1988, Liokumovich-Maximo_2021}. In higher dimensions, this problem remains open.

Interpolating between Ricci and scalar curvature, it is natural to ask whether manifolds with uniformly positive $m$-intermediate curvature have finite Urysohn $(m-1)$-width. Applying Chodosh-Li's slice-and-dice argument \cite{Chodosh-Li_2020}, we show affirmatively the case $m=2$, in dimension up to 5.

\begin{theorem}\label{thm-intro:macroscopic_dim_1}
	Suppose $M$ is a complete simply-connected manifold with dimension $3\leq n\leq 5$ and bi-Ricci curvature $\C_2\geq\lambda>0$. Then $M$ has finite Urysohn 1-width.
\end{theorem}

In larger dimension, on the other hand, there exist counterexamples that are in the same spirit with Theorem \ref{thm-intro:weak-Bonnet-Myers}.

\begin{theorem}\label{thm-intro:countereg_macro_dim_2}
When $n\geq6$, there exists a complete metric on $S^{n-2}\times\RR^2$ with uniformly positive bi-Ricci curvature and infinite Urysohn 1-width.
\end{theorem}

We finally remark that, bi-Ricci curvature lower bound has naturally appeared in the recent work of Chodosh-Li-Minter-Stryker \cite{Chodosh-Li-Minter-Stryker} on the 5-dimensional stable Bernstein theorem. There it is observed that under the Gulliver-Lawson transformation $\tilde g=|x|^{-2}g$, one obtains a metric with $\C_2\geq1$ in the spectral sense, from which one could obtain diameter and volume bounds for warped $\mu$-bubbles.

\vspace{12pt}

This paper is organized as follows. In section \ref{sec:Ricm>=0} we discuss the case of non-negative $\C_m$. This section contains a detailed explanation of the reverse engineering process, and the proof of Theorem \ref{thm-intro:countereg_Ricm}. In section \ref{sec:Ricm>0} we discuss the case of uniformly positive $\C_m$. This section contains an introduction to $\mu$-bubbles, and the proof of Lemma \ref{lemma-intro:weakly_Ric>0}, Theorem \ref{thm-intro:weak-Bonnet-Myers}\,-\,\ref{thm-intro:countereg_macro_dim_2}. In Section \ref{sec:dim6}, we prove the rigidity in dimension 6 (Theorem \ref{thm-intro:rigidity_dim6}).

\vspace{12pt}

\textbf{Notations.} We always use $M$ to denote a smooth oriented manifold, and use $\Sigma$ to denote a hypersurface in $M$. We use $N$ to denote the unit normal vector of $\Sigma$, and $A,H$ denotes the second fundamental form and mean curvature. The sign convention is $A(X,Y)=\metric{\D_X N}{Y}$, $H=\tr_\Sigma A$. These notations apply to all sections except  Section \ref{sec:dim6}, where we adopt the same notations as in the main reference mentioned there.

\vspace{12pt}

\textbf{Acknowledgements.} The author would like to thank Sven Hirsch for inspiring conversations and comments on previous drafts. He also thanks Jianchun Chu and Florian
Johne for useful discussions.

\section{The Case of Nonnegative Intermediate Curvature}\label{sec:Ricm>=0}


\subsection{Reverse engineering of counterexamples}\label{subsec:engi}


In this section, we use the example $T^3\times S^5$ to demonstrate the reverse engineering process mentioned in the introduction part. To better show the idea, we re-write the argument in \cite{Brendle-Hirsch-Johne_2022} using more convenient notations. In what follows, we use subscript to denote the dimension in which an object lives. For example, $u_7$ denotes a function on $M^7$, $N_7$ (resp. $A_6$ and $H_6$) denotes the unit normal vector (resp. second fundamental form and mean curvature) of $M^6\subset M^7$, and $\D_5$ denote the gradient on $M^5$. \\

Assume that $M^8=T^3\times S^5$ has positive tri-Ricci curvature. Ignoring the possible singularity, suppose that we can find a smooth area minimizer $M^7\subset M^8$ in the homology class of $T^2\times S^5$. By the stability inequality for minimal surfaces, there exists $u_7>0$ on $M_7$ such that
\begin{equation}\label{eq-engi:def_u7}
	\Delta_7 u_7 \leq \Big[-|A_7|^2-\Ric_8(N_8,N_8)\Big]u_7.
\end{equation}
Then, let $M^6\subset M^7$ be a minimizer of the functional $M\mapsto\int_M u_7$ in the homology class of $S^1\times S^5$. The stability inequality implies that there exists $u_6>0$ on $M^6$ such that
\begin{equation}\label{eq-engi:def_u6}
	\D_6\cdot\big(u_7\D_6 u_6\big) \leq
	\Big[ -|A_6|^2u_7-\Ric_7(N_7,N_7)u_7+\Delta_7 u_7-\Delta_6 u_7 \Big]u_6.
\end{equation}
Finally, let $M^5\subset M^6$ be a minimizer of $M\mapsto\int_M u_6u_7$ in the homology class of $S^5$. The stability inequality gives
\begin{equation}\label{eq-engi:eq_on_M^5}
	0 \leq \int_{M^5}u_6u_7|\D_5\varphi|^2
	+\int_{M^5}\Big[
		-|A_5|^2u_6u_7-\Ric_6(N_6,N_6)u_6u_7+\Delta_6(u_6u_7)-\Delta_5(u_6u_7)
	\Big]\varphi^2
\end{equation}
for any $\varphi$. We let $\varphi=(u_6u_7)^{-1}$, and then combine (\ref{eq-engi:eq_on_M^5}) with (\ref{eq-engi:def_u6}) (\ref{eq-engi:def_u7}). Note that
\[\Delta_6(u_6u_7) = \D_6\cdot(u_7\D_6u_6)+u_6\Delta_6 u_7+\D_5u_6\cdot\D_5 u_7
	+\frac{\p u_6}{\p N_6}\cdot\frac{\p u_7}{\p N_6},\]
and the Gauss equations
\[\Ric_7(N_7,N_7) = \Ric_8(N_7,N_7)-R_8(N_7,N_8,N_8,N_7)-A_7^2(N_7,N_7),\]
\[\begin{aligned}
	\Ric_6(N_6,N_6) &= \Ric_8(N_6,N_6)-R_8(N_6,N_7,N_7,N_6)-R_8(N_6,N_8,N_8,N_6) \\
	&\qquad +H_6A_6(N_6,N_6)-A_6^2(N_6,N_6)-A_7^2(N_6,N_6) \\
	&\qquad -A_7(N_6,N_6)A_7(N_7,N_7)+A_7(N_6,N_7)^2.
\end{aligned}\]
After rearranging the terms, we obtain
\begin{align}
	0 &\leq -\int_{M^5}\C_3(N_6,N_7,N_8)\varphi \nonumber\\
	&\qquad -\int_{M^5}\Big[(u_6u_7)^{-3}|\D_5(u_6u_7)|^2-(\D_5u_6\cdot\D_5 u_7)(u_6u_7)^{-2}\Big] \nonumber\\
	&\qquad -\int_{M^5}\Big[|A_5|^2-u_6^{-1}\frac{\p u_6}{\p N_6}\cdot u_7^{-1}\frac{\p u_7}{\p N_6}\Big]\varphi \nonumber\\
	&\qquad -\int_{M^5}\Big[|A_6|^2+H_6A_6(N_6,N_6)-A_6^2(N_6,N_6)\Big]\varphi \nonumber\\
	&\qquad -\int_{M^5}\Big[|A_7|^2-A_7^2(N_6,N_6)-A_7^2(N_7,N_7)-A_7(N_6,N_6)A_7(N_7,N_7) \nonumber \\
	&\qquad\qquad\qquad +A_7(N_6,N_7)^2\Big]\varphi \nonumber\\
	&=: -\int_{M^5}\C_3(N_6,N_7,N_8)\varphi-(I+\V_3+\V_2+\V_1). \label{eq-engi:I1-I4}
\end{align}
This inequality is equivalent to Lemma 3.4, 3.10 in \cite{Brendle-Hirsch-Johne_2022}, and there the terms $\V_1,\V_2,\V_3$ are shown to be non-negative in lower dimensions. However, it is not true that $\V_1,\V_2,\V_3\geq0$ in the current dimension. For example, for the term $\V_3$ we have
\[\V_3 \geq \frac{H_5^2}{5}-\frac{\p\log u_6}{\p N_6}\cdot\frac{\p\log u_7}{\p N_6}
= \frac15\Big(\frac{\p\log u_6}{\p N_6}+\frac{\p\log u_7}{\p N_6}\Big)^2
	-\frac{\p\log u_6}{\p N_6}\cdot\frac{\p\log u_7}{\p N_6}
\overset{?}{\geq}0,\]
which is clearly not true by letting $u_6=u_7$ and $A_5$ be a multiple of $g|_{M^5}$. Furthermore, by choosing $A_6,A_7$ of diagonal form with
\begin{equation}\label{eq-engi:countereg_cond2}
	A_6(e_1,e_1)=\cdots=A_6(e_5,e_5)=x,\ \ A_6(N_6,N_6)=-\frac52x,
\end{equation}
and
\begin{equation}\label{eq-engi:countereg_cond3}
	A_7(e_1,e_1)=\cdots=A_7(e_5,e_5)=y,\ \ A_7(N_6,N_6)=A_7(N_7,N_7)=-\frac52y,
\end{equation}
we achieve $\V_2<0$, $\V_1<0$. Among the many available choices, (\ref{eq-engi:countereg_cond2}) (\ref{eq-engi:countereg_cond3}) are made so that the resulting metric has the simplest form. \\

Next, we collect and combine together the information obtained. Consider a metric on $T^3\times S^5$ of the following form:
\[g=f_1^2h+f_2^2dx^2+f_3^2dy^2+f_4^2dz^2,\]
where $f_1,f_2,f_3,f_4$ are functions in $x,y,z$, and $h$ is the standard metric on $S^5$.

The first condition we impose on $g$ is that, any coordinate slice
\[M^5=\big\{x=x_0,y=y_0,z=z_0\big\}\subset M^6=\big\{y=y_0,z=z_0\big\}\subset M^7=\big\{z=z_0\big\}\]
is a slicing of weighted minimal surfaces introduced above. This gives many hints on the choice of $f_i$. First, the mean curvature of $M^7$ is identically zero along the flow $\p/\p z=f_4N_8$ by our condition, hence
\[\Delta_7 f_4\equiv-\big(|A_7^2|+\Ric_8(N_8,N_8)\big)f_4\]
by the variation formula of mean curvature. Compared with (\ref{eq-engi:def_u7}), we have a natural choice $f_4=u_7$. We also know that (\ref{eq-engi:def_u7}) should be equality. Second, the area $\int f_1^5f_2f_3\,dx\,dy$ must be constant in $z$. Hence we let $f_1^5f_2f_3$ be pointwise constant.

Followed by the same analysis on $M^6$, we obtain $f_3=u_6$, and $f_1^5f_2u_7$ is constant (which implies $f_3=f_4=u_7$), and that (\ref{eq-engi:def_u6}) is an equality. By spherical symmetry, (\ref{eq-engi:eq_on_M^5}) is equality for constant $\varphi$. Therefore, all the inequalities (\ref{eq-engi:def_u7}) (\ref{eq-engi:def_u6}) (\ref{eq-engi:eq_on_M^5}) (\ref{eq-engi:I1-I4}) achieve equality, this forces
\[\C_3(N_6,N_7,N_8)=-(\V_1+\V_2+\V_3).\]
Next, we interpret conditions (\ref{eq-engi:countereg_cond2}) (\ref{eq-engi:countereg_cond3}) in terms of $f_1,f_2,f_3,f_4$. Note that
\[A_6(e_i,e_i)=f_3^{-1}\frac{\p\log f_1}{\p y},\ A_6(N_6,N_6)=f_3^{-1}\frac{\p\log f_2}{\p y}\ \xrightarrow{(\ref{eq-engi:countereg_cond2})}\ f_2=f_1^{-5/2},\]
and
\[A_7(e_i,e_i)=f_4^{-1}\frac{\p\log f_1}{\p z},\ A_7(N_7,N_7)=f_4^{-1}\frac{\p\log f_3}{\p z}\ \xrightarrow{(\ref{eq-engi:countereg_cond3})}\ f_3=f_1^{-5/2}.\]
Under these choices we have $\V_1<0,\V_2<0,\V_3<0$. We finally obtain
\[f_1=F(x,y,z),\ f_2=f_3=f_4=F(x,y,z)^{-5/2},\]
which is exactly (\ref{eq-intro:countereg_Ricm}). One could repeat this argument on the general cases $T^m\times S^{n-m}$. A much faster way is to directly consider $g=F^2 h+F^{-s}\sum dx_i^2$. The choice $s=2\frac{n-m}{m-1}$ is determined by cancelling all the second derivatives of $F$ in $\C_m(\p x_1,\cdots,\p x_m)$.

\subsection{Counterexamples in dimension \texorpdfstring{$\geq 7$}{≥7}}

\forcedlinebreak{9pt}

\textbf{Notations}. For convenience, for the rest of this section we let $k=n-m$, so that the underlying manifold is $T^m\times S^k$. Greek indices $\alpha,\beta,\cdots$ are assumed to be in the $S^k$ direction, while Latin indices $i,j,\cdots$ are in the $T^m$ directions. The indices $p,q,\cdots$ can be in either direction. Denote by $\{e_\alpha\}_{\alpha=1}^k$ an orthonormal frame for $g$ tangent to $S^k$, and denote $e_i=\frac{\p x_i}{|\p x_i|}$. Denote the partial derivatives $F_i=\frac{\p F}{\p x_i}$, and $dF^2=\sum_i F_i^2$. Let $s=\frac{2k}{m-1}$. Therefore, the metric under consideration is
\begin{equation}
	g = \epsilon^2F^2h_{\alpha\beta} + F^{-s}\sum_{i=1}^m dx_i^2.
\end{equation}

\begin{proof}[Proof of Theorem \ref{thm-intro:countereg_Ricm}] {\ }
	
We first compute $\C_m\big(\frac{\p x_1}{|\p x_1|},\cdots,\frac{\p x_m}{|\p x_m|}\big)$. The Christoffel symbols of $g$ are:
\[\begin{aligned}
	& \Gamma_{\alpha\beta}^i=-\epsilon^2F^{s+1}F_ih_{\alpha\beta},\quad
	\Gamma_{\alpha i}^\beta=F^{-1}F_i\delta_\alpha^\beta,\quad
	\Gamma_{ij}^\alpha=\Gamma_{i\alpha}^j=0, \\
	& \Gamma_{ij}^i=-\frac s2F^{-1}F_j,\quad 
	(\text{when }i\ne j)\ \Gamma_{ii}^j=\frac s2F^{-1}F_j.
\end{aligned}\]
From this we compute the curvature components (Einstein summation is not used):
\[\begin{aligned}
	& g^{jj}R_{ijj}^i=\frac s2F^{s-1}(F_{ii}+F_{jj})-\frac s2F^{s-2}(F_i^2+F_j^2)-\frac{s^2}4F^{s-2}\sum_{k\ne i,j}F_k^2\quad(i\ne j), \\
	& g^{\alpha\alpha}R_{i\alpha\alpha}^i=-F^{s-1}F_{ii}+\frac s2F^{s-2}\big(\!-F_i^2+\sum_{j\ne i}F_j^2\big), \\
	& g^{\alpha\alpha}R_{i\alpha\alpha}^j=-F^{s-1}F_{ij}-sF^{s-2}F_iF_j\quad(i\ne j), \\
	& R_{ijk}^\alpha=R_{\alpha\beta\gamma}^i=R_{ij\alpha}^\beta=0, \\
	& R_{\alpha\beta\gamma}^\eta=\big(\epsilon^{-2}F^{-2}-F^{s-2}\sum_i F_i^2\big)\big(\delta_\alpha^\eta g_{\beta\gamma}-\delta_\beta^\eta g_{\alpha\gamma}\big).
\end{aligned}\]
Moreover, $R_{ijk}^l$ is always a polynomial combination of $\frac{\p^2 F}{F}$ and $\frac{\p F}F$. Finally we compute
\begin{align}
	\C_m\Big(\frac{\p x_1}{|\p x_1|},\cdots,\frac{\p x_m}{|\p x_m|}\Big)
	&= \sum_{i<j}g^{jj}R_{ijj}^i+\sum_{i,\alpha}g^{\alpha\alpha}R_{i\alpha\alpha}^i \nonumber\\
	&= \Big[\frac{k^2(m-2)}{2(m-1)}-k\Big]F^{s-2}dF^2. \label{eq-engi:Cm_for_torus_factor}
\end{align}
This gives the first statement of Theorem \ref{thm-intro:countereg_Ricm}, and shows that $\C_m$ is nonnegative on any orthonormal frame tangent to $T^m$.
	
Now we suppose $\{X_l\}$ ($1\leq l\leq m$) is an orthonormal frame, not necessarily tangent to $T^m$. We decompose $X_l=X_l^T+X_l^S$, where $X_l^T=\sum_i a_{il}e_i$ are tangent to $T^m$ and $X_l^S=b_{\alpha l}e_\alpha$ are tangent to $S^k$ (hence $\sum_i a_{il}^2+\sum_\alpha b_{\alpha l}^2=1$). We may assume that all the $X_l^T$ factors are nonzero. We compute
\begin{align}
	\C_m(X_1,\cdots,X_m)
	&= \sum_l\Ric(X_l,X_l)-\sum_{l<r}R(X_l,X_r,X_r,X_l) \nonumber\\
	&= \sum_l\Ric(X_l^T,X_l^T)+\sum_l\Ric(X_l^S,X_l^S) \nonumber\\
	&\qquad -\sum_{l<r}R(X_l^T,X_r^T,X_r^T,X_l^T)-\sum_{l<r}R(X_l^S,X_r^S,X_r^S,X_l^S)-J \label{eq-engi:aux1}
\end{align}
where $J$ consists of the cross terms:
\[\begin{aligned}
	J &= \sum_{l<r}\Big[R(X_l^T,X_r^S,X_r^S,X_l^T)+R(X_l^S,X_r^T,X_r^T,X_l^S) \\
	&\qquad +R(X_l^T,X_r^S,X_r^T,X_l^S)+R(X_l^S,X_r^T,X_r^S,X_l^T)\Big].
\end{aligned}\]
From the curvature computations we obtain
\[\Ric(X_l^S,X_l^S) =
	(k-1)\epsilon^{-2}F^{-2}|X_l^S|^2
	+F^s|X_l^S|^2\cdot P\big(\frac{\p F}F,\frac{\p^2 F}{F}\big),\]
where $P(\frac{\p F}F,\frac{\p^2 F}{F})$ is a polynomial combination of $\frac{\p F}F$ and $\frac{\p^2 F}{F}$. Hence
\[\Ric(X_l^S,X_l^S) \geq
	(k-1)\epsilon^{-2}F^{-2}|X_l^S|^2
	-C\big(||\log F||_{C^2}\big)\cdot|X_l^S|^2.\]
Similarly, we estimate
\[|J|\leq C\big(||\log F||_{C^2}\big)\cdot\sum_l|X_l^S|^2.\]
Finally, we have
\[R(X_l^S,X_r^S,X_r^S,X_l^S)\leq\epsilon^{-2}F^{-2}\Big(|X_l^S|^2|X_r^S|^2-\metric{X_l^S}{X_r^S}^2\Big).\]
Hence,
\begin{align}
	\C_m(X_1,\cdots,X_m) &\geq \C_m\big(\frac{X_1^T}{|X_1^T|},\cdots,\frac{X_m^T}{|X_m^T|}\big)-\sum_l\Ric\big(\frac{X_l^T}{|X_l^T|},\frac{X_l^T}{|X_l^T|}\big)|X_l^S|^2 \nonumber\\
	&\qquad +\sum_{l<r}R\big(\frac{X_l^T}{|X_l^T|},\frac{X_r^T}{|X_r^T|},\frac{X_r^T}{|X_r^T|},\frac{X_l^T}{|X_l^T|}\big)\cdot(1-|X_l^T|^2|X_r^T|^2) \nonumber\\
	&\qquad +\epsilon^{-2}F^{-2}\Big[(k-1)\sum_l|X_l^S|^2-\sum_{l<r}\big(|X_l^S|^2|X_r^S|^2-\metric{X_l^S}{X_r^S}^2\big)\Big] \nonumber\\
	&\qquad -C\big(||\log F||_{C^2}\big)\cdot\sum_l|X_l^S|^2. \nonumber
\end{align}
Similarly estimating the second and third term, we have
\begin{align}
	\C_m(X_1,\cdots,X_m) &\geq
		\C_m\big(\frac{X_1^T}{|X_1^T|},\cdots,\frac{X_m^T}{|X_m^T|}\big)
		-C\big(||\log F||_{C^2}\big)\sum_l|X_l^S|^2 \nonumber\\
		&\qquad +\epsilon^{-2}F^{-2}\Big[
			(k-1)\sum_l|X_l^S|^2
			-\sum_{l<r}\big(|X_l^S|^2|X_r^S|^2-\metric{X_l^S}{X_r^S}^2\big)
		\Big]. \nonumber
\end{align}

We note that $\sum_{l=1}^m|X_l^S|^2\leq k$. To see this, we add $k$ vectors $Y_t=Y_t^T+Y_t^S$ ($1\leq t\leq k$) so that $\{X_l\}\cup\{Y_t\}$ forms an orthogonal frame of $T_p M$. We have $\sum_l|X_l^S|^2+\sum_t|Y_t^S|^2=k$ since the coordinates of $\{X_l,Y_t\}$ form an orthogonal matrix. Therefore,
\[\begin{aligned}
	\sum_{l<r}\Big(|X_l^S|^2|X_r^S|^2-\metric{X_l^S}{X_r^S}^2\Big) &\leq
	\frac12\Big(\sum_{l=1}^m|X_l^S|^2\Big)^2-\frac12\sum_{l=1}^m|X_l^S|^4\leq\frac12k\sum_{l=1}^m|X_l^S|^2.
\end{aligned}\]
We finally obtain
\begin{equation}\label{eq-engi:cross_curvatures}
	\begin{aligned}
		\C_m(X_1,\cdots,X_m) &\geq
		\C_m\big(\frac{X_1^T}{|X_1^T|},\cdots,\frac{X_m^T}{|X_m^T|}\big) \\
		&\qquad +\Big[(\frac12k-1)\epsilon^{-2}F^{-2}
			-C\big(||\log F||_{C^2}\big)\Big]\cdot\sum_{l=1}^m|X_l^S|^2.
	\end{aligned}
\end{equation}
Recall that we assume $n(m-2)\geq m^2-2$, this implies $k=n-m\geq\frac{2m-2}{m-2}>2$. Hence the coefficient of $\epsilon^{-2}F^{-2}$ is positive, and the result follows.
\end{proof}

\begin{proof}[Proof of Corollary \ref{cor-intro:perturbation}] {\ }
	
	Let $E$ be the $m$-dimensional Grassmannian bundle over $M=S^k\times T^m$. Thus we may view $\C_m$ as a smooth function on $E$. By the above calculations (\ref{eq-engi:Cm_for_torus_factor}) (\ref{eq-engi:cross_curvatures}), we have $\C_m\geq0$ on $E$, and its minimum on each Grassmannian fiber is attained at $\text{span}(\p_{x_1},\cdots,\p_{x_m})$. Let $Z_0\subset T^m$ be the (finite) set of critical points of $F$. Let $Z\subset E$ be the zero set of $\C_m$. By (\ref{eq-engi:Cm_for_torus_factor}) (\ref{eq-engi:cross_curvatures}), each point $w\in Z$ must have the form
	\[w=\big(y,q;\text{span}(\p_{x_1},\cdots,\p_{x_m})\big),\qquad\text{where $y\in S^k$, $q\in Z_0\subset T^m$.}\]
	For each fixed $w\in Z$, we let $\{x_i\}$ be a translation of the natural coordinates on $T^m$ so that $q=(0,\cdots,0)$, and $\{y_\alpha\}$ be a coordinate chart of $S^k$ near $y$.
	
	The coordinates $\{x_i,y_\alpha\}$ induce a local trivialization of $E$ near $w$. Let $\{z_a\}$ be a coordinate chart of the Grassmannian near $\text{span}(\p_{x_1},\cdots,\p_{x_m})$. By (\ref{eq-engi:Cm_for_torus_factor}) (\ref{eq-engi:cross_curvatures}), the spherical symmetry of the metric, and the vanishing of cross curvature components, the Hessian of $\C_m$ at $w$ has a block diagonal form:
	\[\Hess(\C_m)\big|_w = \begin{pmatrix}
		0 & 0 & 0 \\
		0 & C(m,k)F^{s-2}\Hess^2(F) & 0 \\
		0 & 0 & \geq O(\epsilon^{-2})
	\end{pmatrix},\]
	where the three entries correspond to $y,x,z$, and $\Hess^2(F)_{ij}:=\sum_k\frac{\p^2F}{\p x_i\p x_k}\frac{\p^2F}{\p x_j\p x_k}$. This implies $\C_m\geq C(|x|^2+|z|^2)$ locally near $w$.
	
	Now let $\psi$ be a smooth function on $T^m$ such that $\psi(q_k)=0$, $d\psi(q_k)=0$ and $(\frac{\p^2\psi}{\p x_i\p x_j})_{ij}<0$ at all $q_k$. Consider the perturbation $g_t=(1+t\psi)g$, with $0<t\ll1$. We claim that
	\begin{equation}\label{eq-engi:d/dt(Cm)}
		\frac{d}{dt}\Big[
			\C_m(g_t)(w)\big(\frac{\p x_1}{|\p x_1|_{g_t}},\cdots,\frac{\p x_m}{|\p x_m|_{g_t}}\big)
		\Big]\Big|_{t=0} > 0,\quad\forall\,w\in Z.
	\end{equation}
	To see this, we recall the general formula for variation of metrics: for a family of metrics $g_t$ with $\dot g=\frac{dg}{dt}$, we have $\dot\Gamma_{pq}^r:=\frac{d}{dt}\Gamma_{pq}^r=\frac12g^{rs}(\D_p\dot g_{qs}+\D_q\dot g_{ps}-\D_s\dot g_{pq})$, and $\frac{d}{dt}R_{pqr}^s=\D_p\dot\Gamma_{qr}^s-\D_q\dot\Gamma_{pr}^s$. Applying $\dot g=\psi g$, we obtain $\frac d{dt}R_{pqq}^p=-\frac12(\D^2_{pp}\psi+\D^2_{qq}\psi)$ ($p\ne q$). Note that $\D^2_{ii}\psi<0$ and $\D^2_{\alpha\alpha}\psi=0$ at $w$. Also note that $g_t=g$ at $w$, hence
	\[\frac{d}{dt}\Big[
		\C_m(w)\big(\frac{\p x_1}{|\p x_1|_{g_t}},\cdots,\frac{\p x_m}{|\p x_m|_{g_t}}\big)
	\Big]\Big|_{t=0}
	= \sum_{i<j}g^{jj}\frac d{dt}R_{ijj}^i
		+\sum_{i,\alpha}g^{\alpha\alpha}\frac d{dt}R_{i\alpha\alpha}^i
	>0.\]
	Thus our claim is proved. By the smoothness of all the objects involved, \eqref{eq-engi:d/dt(Cm)} implies that $\frac{d}{dt}\C_m\geq a-b(|x|+|z|)$ in the local coordinate chart defined above, for some $a,b>0$ depending on $F,\psi$. Therefore, there exists a small neighborhood $V_w\ni w$ and a sufficiently small $t_w$ (depending on $F,\psi$) for which
	\[\C_m(g_t) \geq C(|x|^2+|z|^2) + \frac12at - 2bt(|x|+|z|),
	\quad\forall\, t<t_w.\]
	Hence $\C_m(g_t)>0$ in $V_w$ for all sufficiently small $t$. For each $w\in Z$ we obtain such a small neighborhood $V_w$. In the open cover $\bigcup_w V_w\supset Z$, we choose a finite subcover $\{V_{w_i}\}$ and obtain a uniform $t_0\ll1$ such that $\C_m(g_t)>0$ in  $\bigcup V_{w_i}$ when $t\leq t_0$. Since $\C_m(g_0)>0$ on $E\setminus\bigcup V_{w_i}$, by compactness we have $\C_m(g_t)>0$ on $E\setminus\bigcup V_{w_i}$ when $t\leq t_1$, for a sufficiently small $t_1$. Hence $\C_m(g_{\min(t_0,t_1)})>0$ everywhere on $E$.
\end{proof}

\section{The Case of Uniformly Positive Intermediate Curvature}\label{sec:Ricm>0}

\forcedlinebreak{9pt}

\subsection{A brief introduction to \texorpdfstring{$\mu$}{μ}-bubbles}\label{subsec:mu-bubbles}

\forcedlinebreak{9pt}

Let $M$ be a Riemannian manifold. A $\mu$-bubble is a hypersurface $\Sigma\subset M$ with prescribed mean curvature $H=h$, where $h$ is a function on $M$ to be chosen. In variational perspective, we consider the generalized perimeter functional
\begin{equation}\label{eq-mububble:energy_no_weight}
	E_0(\Omega) := |\p\Omega|-\int_\Omega h.
\end{equation}
When $\Omega$ is a critical point of $E_0$, its boundary $\Sigma=\p\Omega$ is a $\mu$-bubble. The second variation formula of $E_0$ at a critical point is computed to be
\begin{equation}\label{eq-mububble:d^2E0}
	\delta^2 E_0(\Omega)(\varphi) = \int_{\Sigma}\Big[
		|\D_\Sigma\varphi|^2 - \big( |A|^2+\Ric(N,N)+\frac{\p h}{\p N} \big)\varphi^2
	\Big].
\end{equation}

A $\mu$-bubble is called \textit{stable} if the second variation (\ref{eq-mububble:d^2E0}) is non-negative for all $\varphi$. To demonstrate how stable $\mu$-bubbles interact with uniformly positive curvatures, here we prove that a stable $\mu$-bubble in a manifold with $\C_2\geq\lambda>0$ satisfies a Bonnet-Myers' theorem, for a properly chosen function $h$.

Let $e$ be a measurable unit vector field on $\Sigma$ such that $\Ric_\Sigma(e,e)=R_0$ almost everywhere. (\ref{eq-mububble:d^2E0}) implies
\[\begin{aligned}
	0 &\leq \int_\Sigma\Big[
		|\D_\Sigma\varphi|^2 + R_0\varphi^2
		- \big( \Ric_\Sigma(e,e)+|A|^2+\Ric_M(N,N)+h_N \big)\varphi^2
	\Big] \\
	&\leq \int_\Sigma\Big[
		|\D_\Sigma\varphi|^2 + R_0\varphi^2 - \C_2(e,N)\varphi^2
		+|\D_M h|\varphi^2 - \big( HA(e,e)-A^2(e,e)+|A|^2 \big)\varphi^2
	\Big].
\end{aligned}\]
When $\dim(M)=n\leq5\Leftrightarrow\dim(\Sigma)\leq4$, we claim that
\[HA(e,e)-A^2(e,e)+|A|^2 \geq C(n)H^2\quad(C(n)>0).\]
\begin{proof}[Proof of the claim]
	Let $\{e_1,\cdots,e_{n-2},e\}$ form an orthonormal frame of $\Sigma$. Then
	\[ HA(e,e)-A^2(e,e)+|A|^2 \geq
		A(e,e)^2 + A(e,e)\sum_{i=1}^{n-2}A(e_i,e_i)
		+ \sum_{i=1}^{n-2} A(e_i,e_i)^2.\]
	Using Young's inequality:
	\[A(e,e)A(e_i,e_i) \geq -\frac45A(e_i,e_i)^2-\frac5{16}A(e,e)^2.\]
	Hence
	\[\begin{aligned}
		HA(e,e)-A^2(e,e)+|A|^2 &\geq
		\min\big( \frac15,1-\frac5{16}(n-2) \big)
		\Big[ A(e,e)^2+\sum_i A(e_i,e_i)^2 \Big] \\
		&\geq C(n)H^2.
	\end{aligned}\]
	(Note: in dimension 6 we can only obtain $HA(e,e)-A^2(e,e)+|A|^2\geq0$.)
\end{proof}

Therefore, we have
\begin{equation}\label{eq-mububble:aux2}
	0\leq\int_\Sigma\Big[
		|\D_\Sigma\varphi|^2+R_0\varphi^2-\frac\lambda2\varphi^2
		-\big(C(n)h^2+\frac\lambda2-|\D_M h|\big)\varphi^2
	\Big],\quad\forall\varphi.
\end{equation}

\begin{lemma}\label{lemma-Ric2>0:bonnet-myers_for_mu_bubble}
	Suppose $\dim(M)\leq 5$, and $M$ satisfies $\C_2\geq\lambda>0$. If the prescribed mean curvature function $h$ satisfies
	\begin{equation}\label{eq-mububble:eq_for_h}
		|\D_M h|\leq C(n)h^2+\frac\lambda2,
	\end{equation}
	then any connected stable $\mu$-bubble satisfies $\lambda_1(-\Delta_\Sigma+R_0)\geq\frac12\lambda$, hence has bounded diameter by Theorem \ref{thm-intro:weak-Bonnet-Myers}. \hfill $\Box$
\end{lemma}

One meaningful choice for $h$ satisfying (\ref{eq-mububble:eq_for_h}) is the following: given a domain $S\subset M$ with smooth boundary, we set
\begin{equation}\label{eq-mububble:choice_of_h}
	h(x)=\sqrt{\lambda/2C(n)}\cot\big[\sqrt{C(n)\lambda/8}\,d(x)\big],
\end{equation}
where $d\in C^\infty$ is a smoothing of $d(-,\p S)$ such that $|\D d|\leq 2$ and $\frac12d(x,\p S)\leq d(x)\leq 2d(x,\p S)$. We note that $h$ is infinite at $\p S$ and away from $S$ by distance $\pi\sqrt{32/C(n)\lambda}=C'(n)\lambda^{-1/2}$. The infinity here plays the role of barrier, so any $\mu$-bubble must lie within distance $C'(n)\lambda^{-1/2}$ from $S$. This observation and Lemma \ref{lemma-Ric2>0:bonnet-myers_for_mu_bubble} are essential for geometric applications of $\mu$-bubbles.

Since our useful choice (\ref{eq-mububble:choice_of_h}) involves infinity, the energy (\ref{eq-mububble:energy_no_weight}) needs to be renormalized in order to be well-defined. Also, in subsequent sections we will need the weighted version of $\mu$-bubbles. Closely following \cite{Chodosh-Li_2020} \cite{Zhu 2020}, we make the definition rigorous as below.

\begin{defn}[weighted $\mu$-bubble]\label{def-mububble:weighted}
	Let $\Omega^+$, $\Omega^-$ be two disjoint domains in $\Sigma$. Suppose $u>0$ is a smooth function. Let $h$ be a smooth prescribed mean curvature function on $\Sigma\setminus(\Omega^-\cup\Omega^+)$, with $h|_{\p\Omega^+}=\infty$, $h|_{\p\Omega^-}=-\infty$. Let $\Omega^0$ be a fixed domain that contains $\Omega^+$ and is disjoint from $\Omega^-$. Consider the following functional acting on all open sets $\Omega$ with $\Omega\Delta\Omega^0\subset\subset\Sigma\setminus(\Omega^-\cup\Omega^+)$:
	\begin{equation}\label{eq-mububble:weighted_energy}
		E_u(\Omega):=\int_{\p\Omega}u\,dl-\int_M(\chi_\Omega-\chi_{\Omega^0})hu\,dA.
	\end{equation}
	A critical point $\Omega$ of $E_u$ is called a $\mu$-bubble with weight $u$.
\end{defn}

We will usually call $\p\Omega$ a stable $\mu$-bubble as well, as no ambiguity occurs. It is shown in \cite{Chodosh-Li_2020} \cite{Zhu 2020} that a global minimizer of (\ref{eq-mububble:weighted_energy}) always exists. Since $h$ is infinite on $\p\Omega^\pm$, a $\mu$-bubble always lies between $\Omega^+$ and $\Omega^-$. By geometric measure theory, a $\mu$-bubble is always a $C^{2,\alpha}$ hypersurface. The first variation of (\ref{eq-mububble:weighted_energy}) gives
\begin{equation}\label{eq-appA: solve1}
	H = h-u^{-1}u_N,
\end{equation}
and the second variation of (\ref{eq-mububble:weighted_energy}) at a critical point gives
\begin{equation}\label{eq-mububble:2nd_variation}
	\begin{aligned}
		\delta^2 E_u(\varphi) &= \int_\Sigma\Big[
			u|\D_\Sigma\varphi|^2-\big(|A|+\Ric(N,N)+h_N\big)u\varphi^2 \\
			&\qquad\qquad +(\Delta u-\Delta_\Sigma u)\varphi^2-hu_N\varphi^2
		\Big].
	\end{aligned}
\end{equation}
\begin{proof}[Proof of (\ref{eq-mububble:2nd_variation})]
	The first variation of (\ref{eq-mububble:weighted_energy}) is computed to be
	\[\delta E_u(\varphi) = \int_\Sigma\big( Hu+u_N-hu \big)\varphi\]
	Since $\Sigma$ is a critical point of $E_u$, the second variation contains the following terms:
	\begin{equation}\label{eq-mububble:aux1}
		\delta^2 E_u(\varphi) = \int_\Sigma\Big[
			\frac{dH}{dt}u + Hu_N\varphi + \frac{du_N}{dt}
			-h_N u\varphi - hu_N\varphi
		\Big]\varphi
	\end{equation}
	Finally, it is an elementary fact that
	\[\frac{dH}{dt} = -\Delta\varphi-\big(|A|^2+\Ric(N,N)\big)\varphi,\quad
	\frac{du_N}{dt} = \D^2 u(N,N)\varphi-\metric{\D_\Sigma u}{\D_\Sigma\varphi}.\]
	Combining these with (\ref{eq-mububble:aux1}) and using integration by part, we obtain (\ref{eq-mububble:2nd_variation}).
\end{proof}

\textbf{Notes.} $\mu$-bubbles were first considered by Gromov \cite{Gromov_2018} \cite{Gromov_Lectures} to prove metric inequalities on scalar curvature. Recently, this technique has found many applications in problems involving scalar curvature. J. Zhu \cite{Zhu 2020} utilized $\mu$-bubbles to obtain several geometric inequalities. Lesourd-Unger-Yau \cite{Lesourd-Unger-Yau_2021} applied the $\mu$-bubble technique to prove a positive mass theorem with arbitrary ends. Chodosh-Li \cite{Chodosh-Li_2020} and Gromov \cite{Gromov_2020 no PSC} recently proved that aspherical 4- and 5-manifolds do not admit metrics with positive scalar curvature.

\subsection{Results on diameter bounds}\label{subsec:weak_bonnet_myers}

\forcedlinebreak{9pt}

\begin{proof}[Proof of Lemma \ref{lemma-intro:weakly_Ric>0}]
	Let $\Sigma\subset M$ as in the statement of lemma. The stability inequality gives
	\begin{equation}\label{eq-Ric2>0:2nd_variation}
		0 \leq \int_\Sigma\Big[
			|\D\varphi|^2-\big(|A|^2+\Ric_M(N,N)\big)\varphi^2
		\Big], \quad\forall\varphi\in C^\infty(\Sigma),
	\end{equation}
	where we use $N$ to denote the unit normal vector. We make a measurable choice of unit vector field $e$ on $\Sigma$, such that $\Ric_\Sigma(e,e)=R_0$ at every point. By Gauss equation $\Ric_\Sigma(e,e)=\Ric_M(e,e)-R_M(N,e,e,N)-A^2(e,e)$, from (\ref{eq-Ric2>0:2nd_variation}) we obtain
	\begin{align*}
		0 &\leq \int_\Sigma\Big[
			|\D\varphi|^2+R_0\varphi^2-\big(\Ric_\Sigma(e,e)+|A|^2+\Ric_M(N,N)\big)\varphi^2
		\Big] \\
		&= \int_\Sigma\Big[
			|\D\varphi|^2+R_0\varphi^2-\C_{2,M}(e,N)\varphi^2-\big(|A|^2-A^2(e,e)\big)\varphi^2
		\Big] \\
		&\leq \int_\Sigma\Big[
			|\D\varphi|^2+R_0\varphi^2-\lambda\varphi^2
		\Big], \quad\forall\varphi\in C^\infty(\Sigma).\qedhere
	\end{align*}
\end{proof}

\forcedlinebreak{6pt}

\begin{proof}[Proof of Theorem \ref{thm-intro:weak-Bonnet-Myers}] {\ }
	
	\forcedlinebreak{6pt}
	
	\textit{Part 1: proof of the diameter bound.}
	
	Given the eigenvalue condition in the statement of the theorem, we let $u$ satisfy
	\[\Delta u \leq \beta R_0u-\lambda u,\]
	Consider the weight function $v=u^{1/\beta}$, which satisfies the inequality
	\[\Delta v \leq R_0v-\lambda\beta^{-1}v+(1-\beta)v^{-1}|\D v|^2.\]
	
	Let $h$ be a prescribed mean curvature function to be determined. Let $S\subset\Sigma$ be a stable $\mu$-bubble for the energy functional (\ref{eq-mububble:weighted_energy}) with weight $v$. Then $S$ satisfies the mean curvature condition (denote $v_N=\frac{\p v}{\p N}$ for short)
	\begin{equation}\label{eq-Ric2>0:solve1}
		H=h-v^{-1}v_N,
	\end{equation}
	as well as the stability inequality
	\[0 \leq \int_S\Big[
		v|\D_S\varphi|^2 - (|A|^2+\Ric_\Sigma(N,N)+h_N)v\varphi^2
		+(\Delta_\Sigma v-\Delta_S v)\varphi^2-hv_N\varphi^2
	\Big],\ \ \forall\varphi.\]
	We will derive a contradiction from the stability inequality. Replacing $\varphi=v^{-1}$ we obtain
	\begin{align}
		0 &\leq \int_S\Big[
			-v^{-3}|\D_S v|^2 - \big(\frac{(h-v^{-1}v_N)^2}{n-1}+R_0+h_N\big)v^{-1} \nonumber\\
			&\qquad +R_0v^{-1}-\lambda\beta^{-1}v^{-1}
			+(1-\beta)v^{-3}\big(|\D_S v|^2+v_N^2\big)-hv^{-2}v_N
		\Big] \label{eq-Ric2>0:stability_ineq}\\
		&= \int_S\Big[
			-\beta v^{-3}|\D_S v|^2 - \big(\frac1{n-1}h^2+\lambda\beta^{-1}+h_N\big)v^{-1} \nonumber\\
			&\qquad +(1-\beta-\frac1{n-1})v^{-3}v_N^2
			-\frac{n-3}{n-1}hv^{-2}v_N
		\Big] \label{eq-Ric2>0:stability_ineq2}
	\end{align}
	Under the main condition (\ref{eq-intro:dim-beta_ineq}), the coefficient of $v^{-3}v_N^2$ is negative. When $n\ne 3$, we further apply Young's inequality to the last term:
	\begin{equation}\label{eq-Ric2>0:solve2}
		-\frac{n-3}{n-1}hv^{-2}v_N \leq
			(\frac1{n-1}+\beta-1)v^{-3}v_N^2
			+\frac14\big(\frac{n-3}{n-1})^2\frac{h^2v^{-1}}{\frac1{n-1}+\beta-1}.
	\end{equation}
	Therefore, a contradiction is obtained if $h$ satisfies
	\begin{equation}\label{eq-Ric2>0:solve3}
		|\D_\Sigma h|<\Big[
		\frac1{n-1} - \frac14\big(\frac{n-3}{n-1})^2\frac{1}{\frac1{n-1}+\beta-1}
		\Big]h^2+\lambda\beta^{-1}
		=:C_2h^2+\lambda\beta^{-1}.
	\end{equation}
	Note that (\ref{eq-intro:dim-beta_ineq}) implies $C_2>0$.
	Let $p\in M$, and $\Omega^+$ be a small geodesic ball near $p$. Set
	\[h(x) = \sqrt{\lambda\beta^{-1}/C_2}
		\cot\Big[ \sqrt{C_2\lambda\beta^{-1}/4}\,d(x) \Big],\]
	where $d(x)$ is a smoothing of $d(x,\p\Omega^+)$ such that $|\D d|\leq2$, $\frac12d(x,\p\Omega^+)\leq d(x)\leq 2d(x,\p\Omega^+)$. Then $h$ satisfies (\ref{eq-Ric2>0:solve3}) under $C_2>0$. If $\diam(\Sigma)>4\pi/\sqrt{C_2\lambda\beta^{-1}}$, then $\Omega^-=\{q\in\Sigma: d(q,p)>2\pi/\sqrt{C_2\lambda\beta^{-1}}\}$ is non-empty for some choice of $p$. Then the $\mu$-bubble problem (\ref{eq-mububble:weighted_energy}) is well-defined, and we obtain a contradiction from the above argument. This shows that $\diam(\Sigma)\leq4\pi/\sqrt{C_2\lambda\beta^{-1}}$ must hold.
	
	\begin{remark}\label{rmk-Ric2>0:top_obstruction}
		We can show that $S^{n-1}\times S^1$ ($n\leq7$) does not admit metrics with $\lambda_1(-\Delta+\beta R_0)>0$, when $\beta\geq1-\frac1{n-1}$. Otherwise, we find an minimizer of the weighted area $\Sigma\mapsto\int_\Sigma v$ in the homology class of $S^{n-1}$, and obtain contradiction from (\ref{eq-Ric2>0:stability_ineq2}) with $h=0$.
	\end{remark}

	\vspace{12pt}
	
	\textit{Part 2: construction of counterexamples.}
	
	The hyperbolic space has first Laplacian eigenvalue $\frac14(n-1)^2$, therefore is a counterexample for $\beta<\frac{n-1}4$. For the critical case $\beta=\frac14(n-1)$, we use reverse-engineering to find warped product metrics $(S^{n-1}\times\RR,g)$ with
	\[g = dr^2+\epsilon^2f(r)^2\bar g,\quad
		r\in(-\infty,\infty),\ \ 
		\bar g=\text{standard metric on }S^{n-1},\]
	that satisfies (\ref{eq-intro:weak_Ricci_lower_bound}). This process also produce interesting examples with finite volume, in contrast with the hyperbolic plane. Below we still assume in general that $\beta\leq\frac{n-1}4$ when $n\ne3$ and $\beta<\frac12$ when $n=3$.
	
	The metric $g$ above, together with the first eigenfunction $u=u(r)$, will be constructed to satisfy
	\begin{equation}\label{eq-Ric2>0:countereg_cond}
		\Delta u=\beta\Ric(\p_r,\p_r)u-\lambda u
		\ \ \text{and}\ \ 
		\Ric(\p_r,\p_r)\leq\Ric(e,e),\ \forall\,|e|=1.
	\end{equation}
	Thus $g$ is a non-compact counterexample. It is not hard to compute
	\begin{equation}\label{eq-Ric2>0:curvatures}
		\Ric(\p_r,\p_r) = -(n-1)\frac{f''}f,\quad
		\Ric(e',e') = (n-2)\epsilon^{-2}f^{-2} - \frac{(n-2)f'^2+ff''}{f^2},
	\end{equation}
	where $e'$ is any unit vector tangent to $S^{n-1}$. In light of the $\mu$-bubble proof above, it is convenient to switch to the function
	\begin{equation}\label{eq-Ric2:uv_rel}
		v=u^{1/\beta}.
	\end{equation}
	The first part of (\ref{eq-Ric2>0:countereg_cond}) is thus equivalent to
	\begin{equation}\label{eq-Ric2:eq_for_v}
		v'' + (n-1)\frac{f'}fv' + (n-1)\frac{f''}fv + \lambda\beta^{-1}v + (\beta-1)v^{-1}(v')^2 = 0.
	\end{equation}
	To start with, we require that all the constant $r$ spheres are $\mu$-bubbles, similar to subsection \ref{subsec:engi}. Hence we impose the relation
	\begin{equation}\label{eq-Ric2>0:eq1_for_v}
		(n-1)\frac{f'}f = h-\frac{v'}v.
	\end{equation}
	
	Next we collect the inequalities used in the proof and solve the equality case. The case $n=3$ is special since the inequalities (\ref{eq-Ric2>0:solve2}) (\ref{eq-Ric2>0:solve3}) above are not used in the proof. In this case, following directly from (\ref{eq-Ric2>0:stability_ineq2}) we impose the following constraint equation:
	\begin{equation}\label{eq-Ric2>0:eq_dim3}
		(\frac12-\beta)v^{-2}(v')^2 = \frac12h^2+\lambda\beta^{-1}+h'.
	\end{equation}
	It is elementary to verify that any functions $f,v,h$ that satisfy (\ref{eq-Ric2>0:eq1_for_v}) and (\ref{eq-Ric2>0:eq_dim3}) will also satisfy (\ref{eq-Ric2:eq_for_v}). The choice is not unique, among which we consider
	\begin{equation}\label{eq-Ric2>0:sol_dim3}
		\left\{\begin{aligned}
			h(r) &= -\lambda\beta^{-1}r, \\
			v(r) &= \exp\big( \frac{\lambda\beta^{-1} r^2}{2\sqrt{1-2\beta}} \big), \\
			f(r) &= \exp\Big[ -\frac14(\frac1{\sqrt{1-2\beta}}+1) \lambda\beta^{-1}r^2 \Big].
		\end{aligned}\right.
	\end{equation}
	For this choice $f$ decays super-exponentially at infinity. By the computations (\ref{eq-Ric2>0:curvatures}), the second part of (\ref{eq-Ric2>0:countereg_cond}) holds for sufficiently small $\epsilon$. This finishes the construction for $n=3$.
	
	Next we assume $n\ne 3$. Instead of considering the equality case in (\ref{eq-Ric2>0:solve3}), it turns out to be more convenient to impose the relation
	\begin{equation}\label{eq-Ric2>0:eq2_for_v}
		v'=C_3hv,\quad\text{where}\ \,C_3=\frac2{3-n}.
	\end{equation}
	In light of (\ref{eq-Ric2>0:stability_ineq2}) it is then natural to consider the equation
	\begin{align}
		h' &= -\lambda\beta^{-1}+\Big[-\frac1{n-1}-\frac{n-3}{n-1}C_3+(1-\beta-\frac1{n-1})C_3^2\Big]h^2 \nonumber \\
		&= -\lambda\beta^{-1}+\frac{n-1-4\beta}{(n-3)^2}h^2. \label{eq-Ric2>0:eq3_for_v}
	\end{align}
	It can be verified that the solution to (\ref{eq-Ric2>0:eq1_for_v}) (\ref{eq-Ric2>0:eq2_for_v}) (\ref{eq-Ric2>0:eq3_for_v}) satisfies (\ref{eq-Ric2:eq_for_v}).
	In the limit case $\beta=\frac{n-1}4$, $n\ne 3$, we obtain an exact solution
	\begin{equation}\label{eq-Ric2>0:sol_limit}
		h(r) = -\lambda\beta^{-1}r,\quad
		v(r) = \exp\big(\frac{\lambda\beta^{-1}r^2}{n-3}\big),\quad
		f(r) = \exp\big(-\frac{\lambda\beta^{-1}r^2}{2(n-3)}\big).
	\end{equation}
	When $n\geq4$ the minimal Ricci eigenvalue is achieved by $\Ric(\p_r,\p_r)$ for sufficiently small $\epsilon$, for the same reason as above. When $n=2$ all the Ricci curvature is equal to the Gauss curvature. Hence the second part of (\ref{eq-Ric2>0:countereg_cond}) holds for small $\epsilon$. When $\beta<\frac{n-1}4$, we obtain a solution
	\begin{equation}\label{eq-Ric2>0:sol_non_limit}
		\left\{\begin{aligned}
			h(r) &= -\sqrt{\lambda\beta^{-1}/C_4}
			\tanh\Big[ \sqrt{C_4\lambda\beta^{-1}}r \Big] \\
			v(r) &= \cosh\Big[ \sqrt{C_4\lambda\beta^{-1}}r \Big]^{-C_3/C_4} \\
			f(r) &= \cosh\Big[ \sqrt{C_4\lambda\beta^{-1}}r \Big]^{(C_3-1)/C_4(n-1)}
		\end{aligned}\right.
	\end{equation}
	where $C_4:=\frac{n-1-4\beta}{(n-3)^2}>0$. Since $C_3<0$ when $n\geq4$, the second part of (\ref{eq-Ric2>0:countereg_cond}) holds for small $\epsilon$. To summarize, the solutions (\ref{eq-Ric2:uv_rel}) (\ref{eq-Ric2>0:sol_dim3}) (\ref{eq-Ric2>0:sol_limit}) (\ref{eq-Ric2>0:sol_non_limit}) satisfy (\ref{eq-Ric2>0:countereg_cond}), and are easily verified to have finite volume when $n\geq3$. This completes the construction.
\end{proof}

\begin{remark}
	In the non-critical case $\beta<\frac14(n-1)$ when $n\geq4$, one is allowed to vary the parameter $C_3$ as long as $C_4=-\frac1{n-1}-\frac{n-3}{n-1}C_3+(1-\beta-\frac1{n-1})C_3^2$ remains nonnegative. When $\beta\geq 1-\frac1{n-1}$ all valid choices satisfy $C_3\leq\frac1{2-n}$, and the resulting metric $g$ has finite volume. When $\beta<1-\frac1{n-1}$ there exist choices $C_3>1$ that correspond to expanding end at infinity.
\end{remark}

\subsection{Results on Macroscopic Dimension}

\forcedlinebreak{9pt}

\begin{proof}[Proof of Theorem \ref{thm-intro:macroscopic_dim_1}] {\ }
	
	The proof is an application of Chodosh and Li's slice-and-dice argument \cite{Chodosh-Li_2020}.
	
	Let $\Omega_0\subset M$ be a small geodesic ball. Consider the minimizer of the energy functional (\ref{eq-mububble:weighted_energy}) with unit weight $u=1$ and mean curvature function $h$ given by (\ref{eq-mububble:choice_of_h}) with $S=\Omega_0$. We obtain a stable $\mu$-bubble $\Omega_1$ with $\Omega_0\subset\Omega_1\subset N(\Omega_0,C(n)\lambda^{-1/2})$, where $N(\Omega,c)$ denotes the $c$-distance neighborhood of $\Omega$. Replacing $\Omega_1$ by its connected component containing $\Omega_0$, we may assume $\Omega_1$ is connected. Then consider the same $\mu$-bubble problem but with the choice of domain $S=N(\Omega_1,1)$. We obtain another stable $\mu$-bubble $\Omega_2$ with $N(\Omega_1,1)\subset\Omega_2\subset N(\Omega_1,1+C(n)\lambda^{-1/2})$. Repeating this process, we obtain an exhausting chain of domains $\Omega_0\subset\Omega_1\subset\Omega_2\subset\cdots$, with each $\Omega_i$ connected. Denote $\Sigma_i=\p\Omega_i$. Let $\Sigma_{ik}$ be any connected component of $\Sigma_i$. By simply-connectedness of $M$, $\Sigma_{ik}$ is separating, i.e. $M\setminus\Sigma_{ik}$ has two connected components. One component of $M\setminus\Sigma_{ik}$ contains $\Omega_i$ while the other one does not.
	
	Let $M\setminus(\bigcup_i\Sigma_i)=\bigcup_{ij}U_{ij}$ be a decomposition into connected components, with $U_{ij}\subset\Omega_{i+1}\setminus\Omega_i$. Note that $\p U_{ij}$ contains exactly one connected component $\Sigma_{ij}$ of $\Sigma_i$. (The adjacency pattern of $U_{ij}$ is thus a tree.)
	
	By construction we have $U_{ij}\subset N(\Sigma_i,1+C(n)\lambda^{-1/2})$. By the separating property of $\Sigma_{ij}$, we actually have $U_{ij}\subset N(\Sigma_{ij},1+C(n)\lambda^{-1/2})$. By Lemma \ref{lemma-Ric2>0:bonnet-myers_for_mu_bubble} we have $\diam(\Sigma_{ij})\leq C(n)\lambda^{-1/2}$. Hence $\diam(U_{ij})$ is uniformly bounded. We thicken the set $U_{ij}$ by adding a small distance neighborhood of $\Sigma_{ij}$: let $\tilde U_{ij}=U_{ij}\cup N(\Sigma_{ij},\epsilon_{ij})$. Thus $\{\tilde U_{ij}\}$ form an open cover. When $\epsilon_{ij}\ll1$ and $\epsilon_{ij}\ll\min_k d(\Sigma_{ij},\Sigma_{ik})$, there is no overlap among three sets. Let $X$ be the nerve of $\{\tilde U_{ij}\}$, thus $X$ is a 1-dimensional complex.
	
	Finally, let $\{f_{ij}\}$ be a partition of unity subordinate to $\{\tilde U_{ij}\}$. By linear interpolation with weight $f_{ij}$, we obtain a continuous map $\Phi: M\to X$. For each point $p$, the preimage $\Phi^{-1}(p)$ is contained in a single open set in $\{\tilde U_{ij}\}$, hence has uniformly bounded diameter.
\end{proof}

\begin{proof}[Proof of Theorem \ref{thm-intro:countereg_macro_dim_2}] {\ }

	To avoid notational conflicts, we use $k$ to denote the dimension in this proof. Consider $M=\RR^2\times S^{k-2}$ ($k\geq6$) with the metric
	\[g = u(r)^2dt^2 + dr^2 + \epsilon^2f(r)^2\bar g,\]
	where $\bar g$ denotes the round metric on $S^{k-2}$. For the functions $u(r)$ and $f(r)$, we choose (\ref{eq-Ric2:uv_rel}) (\ref{eq-Ric2>0:sol_limit}) for $k=6$ and (\ref{eq-Ric2:uv_rel}) (\ref{eq-Ric2>0:sol_non_limit}) for $k\geq7$. Both choices are applied to the parameters $\beta=1$, $n=k-1$. Note that $u\geq1$ in these choices, which implies that $g$ has infinite Urysohn 1-width: for any point $p\in S^{k-2}$, the map $(\RR^2,g_{Euc})\to\RR^2\times\{p\}\subset(M,g)$ is expanding. If $g$ has finite Urysohn 1-width, then so does $\RR^2$. But $\RR^2$ is known to have infinite Urysohn 1-width.
	
	It remains to confirm that $g$ has uniformly positive bi-Ricci curvature. We compute the curvature components of $g$:
	\begin{equation}\label{eq-Ric2>0:aux5}
		\begin{aligned}
			& \sec(u^{-1}\p_t,\p_r) = -\frac{u''}u,\quad
			\sec(\p_r,u^{-1}\p_t) = -\frac{u''}u,\quad
			\sec(e,e') = \epsilon^{-2}f^{-2}-\frac{(f')^2}{f^2}, \\
			&\qquad\qquad\qquad\quad
			\sec(\p_r,e) = -\frac{f''}f,\quad
			\sec(u^{-1}\p_t,e) = -\frac{u'f'}{uf}.
		\end{aligned}
	\end{equation}
	where $e\perp e'$ are any two unit vectors tangent to $S^{k-2}$. All the cross curvature terms $R_{pqrs}$ ($q\ne r$) vanish. As expected from the construction (see equation (\ref{eq-Ric2:eq_for_v})), we have
	\[\C_2(\p_r,u^{-1}\p_t)
		= -\frac{u''}u
		-(k-2)\frac{f''}f
		-(k-2)\frac{u'f'}{uf}\equiv\lambda.\]
	To control $\C_2(X_1,X_2)$ where $X_1,X_2$ are not tangent to $\p_r,\p_t$, we repeat the argument in (\ref{eq-engi:cross_curvatures}) and above. The estimates of curvature components are replaced by (\ref{eq-Ric2>0:aux5}). Note that $\epsilon^{-2}f^{-2}$ has at least exponential growth at $r\to\infty$, while all the other terms in (\ref{eq-Ric2>0:aux5}) have at most polynomial growth. Following the argument from (\ref{eq-engi:aux1}) to (\ref{eq-engi:cross_curvatures}), we confirm $\C_2\geq\lambda$ when $\epsilon$ is sufficiently small.
\end{proof}

\section{Rigidity in Dimension 6}\label{sec:dim6}

\begin{proof}[Proof of Theorem \ref{thm-intro:rigidity_dim6}] {\ }

	In the proof of Theorem \ref{thm-intro:Chu-Kwong-Lee} in \cite{Chu-Kwong-Lee_2022} (which is labeled Theorem 1.1 there), all the statements only assume $n(m-2)<m^2-2$ except for Proposition 3.1 where a narrower inequality $n(m-1)<m(m+1)$ is assumed. The former inequality holds whenever $n\leq6$, while the latter requires $n\leq5$. The presence of inequality $n(m-1)<m(m+1)$ is due to the usage of $\mu$-bubbles in the proof of equation (3.8) in \cite{Chu-Kwong-Lee_2022} (see (3.11) and its context). As we have previously noticed, $\mu$-bubbles narrow the constraints on the dimension.
	
	Here we find another proof of the equation \cite[(3.8)]{Chu-Kwong-Lee_2022} without using $\mu$-bubbles, which is applicable in dimension $n=6$. The remaining parts of the proof in \cite{Chu-Kwong-Lee_2022} are left unchanged. As a result, this extends Theorem \ref{thm-intro:Chu-Kwong-Lee} to dimension 6 and proves Theorem \ref{thm-intro:rigidity_dim6}. Our argument is inspired by the scalar curvature rigidity theorem of Bray-Brendle-Neves \cite{Bray-Brendle-Neves_2010} (see also \cite[Section 2.4]{Lee_GeomRela} and references therein).
	
	For consistency, in this section we adopt the same notations as in \cite{Chu-Kwong-Lee_2022}. We briefly summarize the setups for the equation \cite[(3.8)]{Chu-Kwong-Lee_2022}. Let $M=T^m\times M^{n-m}$ be a manifold with $\C_m\geq0$. We inductively construct a chain of hypersurfaces $M=\Sigma_0\supset\Sigma_1\supset\cdots\supset\Sigma_m$ and positive weight functions $\rho_i\in C^\infty(\Sigma_0)$ that satisfies the conditions in \cite[Definition 2.1]{Chu-Kwong-Lee_2022}. The precise construction is as follows. We start with setting $\rho_0=1$ and $\Sigma_0=M$. For each $1\leq k\leq m$, we inductively define $\Sigma_k$ as an area minimizer of the functional
	\[\H_{\rho_{k-1}}^{n-k}(\Sigma)=\int_\Sigma\rho_{k-1}\]
	in a certain homology class in $\Sigma_{k-1}$, then set $v_k$ to be the first eigenfunction of the associated stability operator.	Finally, inductively set the weight function $\rho_k=v_k\rho_{k-1}$.
	
	By the first variation formula, the generalized mean curvature
	\begin{equation}\label{eq-dim6:gen_H}
		H_{\rho_{k-1}}(\Sigma_k):=H_{\Sigma_k}+\metric{\D_{\Sigma_{k-1}}\log\rho_{k-1}}{\nu_k}
	\end{equation}
	vanishes on each $\Sigma_k$, where $\nu_k$ is the normal vector of $\Sigma_k\subset\Sigma_{k-1}$. By analyzing the stability inequalities for the generalized area functionals $\H_{\rho_{k-1}}^{n-k}$, a list of rigidity statements are obtained on $\Sigma_k$, see \cite[Proposition 2.1]{Chu-Kwong-Lee_2022}. Next, a foliation of hypersurfaces $\{\Sigma_{m,t}\subset\Sigma_{m-1}\}_{0\leq t\leq\epsilon}$, with $\Sigma_{m,0}=\Sigma_m$, is constructed such that the generalized mean curvature $H_{\rho_{m-1}}(\Sigma_{m,t})$ is constant on each $\Sigma_{m,t}$. Let $\varphi_t\nu_m$ be the corresponding variational vector field, and denote by $H_\rho(t)$ the generalized mean curvature of $\Sigma_{m,t}$. From the construction in \cite[Lemma 3.1]{Chu-Kwong-Lee_2022}, we infer that $\varphi_0=1$, hence $\varphi_t>0$ for sufficiently small $t$. Equation \cite[(3.8)]{Chu-Kwong-Lee_2022} states that
	\begin{equation}\label{eq-dim6:H<0}
		H_\rho(t)\leq0
	\end{equation}
	for all $t$ where the foliation is defined. Based on \eqref{eq-dim6:H<0}, the weighted volume of $\Sigma_{m,t}$ is non-increasing. By the minimality of $\Sigma_{m,0}=\Sigma_m$, all the slices $\Sigma_{m,t}$ are minimizing as well, and thus the rigidity statements in \cite[Proposition 2.1]{Chu-Kwong-Lee_2022} propagate along the foliation. From this one obtains \cite[Proposition 3.1]{Chu-Kwong-Lee_2022}, which finally implies \cite[Theorem 1.1]{Chu-Kwong-Lee_2022} using a continuation argument.
	
	The new proof of (\ref{eq-dim6:H<0}) is obtained by directly computing $dH_\rho(t)/dt$. By the classical variation formula for mean curvature and the definition \eqref{eq-dim6:gen_H}, we have
	\begin{equation}\label{eq-dim6:aux1}
		\begin{aligned}
			\frac{dH_\rho(t)}{dt} &=
				-\Delta_{\Sigma_m}\varphi_t-\big(|h_{\Sigma_m}|^2+\Ric_{\Sigma_{m-1}}(\nu_m,\nu_m)\big)\varphi_t \\
				&\qquad\quad +\varphi_t\D^2_{\Sigma_{m-1}}\log\rho_{m-1}(\nu_m,\nu_m)
				-\metric{\D_{\Sigma_m}\varphi_t}{\D_{\Sigma_m}\log\rho_{m-1}}
		\end{aligned}
	\end{equation}
	where we write $\Sigma_m=\Sigma_{m,t}$ for short. Multiplying both sides by $\varphi_t^{-1}$ and integrating over $\Sigma_{m,t}$, we obtain
	\[\begin{aligned}
		\frac{dH_\rho(t)}{dt}\int_{\Sigma_{m,t}}\varphi_t^{-1}
		&= \int_{\Sigma_{m,t}}\Big[
			-\varphi_t^{-1}\Delta_{\Sigma_m}\varphi_t
			-|h_{\Sigma_m}|^2-\Ric_{\Sigma_{m-1}}(\nu_m,\nu_m) \\
			&\qquad
			+\Delta_{\Sigma_{m-1}}\log\rho_{m-1}-\Delta_{\Sigma_m}\log\rho_{m-1}
			-H_{\Sigma_m}\metric{\D_{\Sigma_{m-1}}\log\rho_{m-1}}{\nu_m} \\
			&\qquad
			-\metric{\D_{\Sigma_m}\log\varphi_t}{\D_{\Sigma_m}\log\rho_{m-1}}
		\Big]
	\end{aligned}\]
	Set $\tilde\rho_m=\varphi_t\rho_{m-1}$. From integration by part we obtain
	\[\int_{\Sigma_{m,t}}\Big[
		-\varphi_t^{-1}\Delta_{\Sigma_m}\varphi_t
		-\metric{\D_{\Sigma_m}\log\varphi_t}{\D_{\Sigma_m}\log\rho_{m-1}}
	\Big] = -\int_{\Sigma_{m,t}}
		\metric{\D_{\Sigma_m}\log\varphi_t}{\D_{\Sigma_m}\log\tilde\rho_m}\]
	Therefore,
	\[\begin{aligned}
		\frac{dH_\rho(t)}{dt}\int_{\Sigma_{m,t}}\varphi_t^{-1}
		&= \int_{\Sigma_{m,t}}\Big[
			-|h_{\Sigma_m}|^2-\Ric_{\Sigma_{m-1}}(\nu_m,\nu_m)
			-\metric{\D_{\Sigma_m}\log\varphi_t}{\D_{\Sigma_m}\log\tilde\rho_m} \\
			&\qquad +\Delta_{\Sigma_{m-1}}\log\rho_{m-1}
			-H_{\Sigma_m}\metric{\D_{\Sigma_{m-1}}\log\rho_{m-1}}{\nu_m}
		\Big]
	\end{aligned}\]
	We need to show that the right hand side is non-positive. This is done by modifying the main argument in \cite{Brendle-Hirsch-Johne_2022} (the necessity of modification is due to the fact that the bottom slice $\Sigma_{m,t}$ here has nonzero generalized mean curvature). The slicing identities \cite[Lemma 3.1, 3.2]{Brendle-Hirsch-Johne_2022} are applicable to our situation for $1\leq k\leq m-1$ since they do not involve the bottom slice $\Sigma_{m,t}$. Applying these identities to deform the term $\Delta_{\Sigma_{m-1}}\log\rho_{m-1}$ (see the proof of \cite[Lemma 3.4]{Brendle-Hirsch-Johne_2022}), we obtain
	\begin{align}
		\frac{dH_\rho(t)}{dt}\int_{\Sigma_{m,t}}\varphi_t^{-1}
		&\leq \int_{\Sigma_{m,t}}\Big[
			-|h_{\Sigma_m}|^2-\Ric_{\Sigma_{m-1}}(\nu_m,\nu_m)
			-\metric{\D_{\Sigma_m}\log\varphi_t}{\D_{\Sigma_m}\log\tilde\rho_m} \nonumber\\
			&\qquad -H_{\Sigma_{m,t}}\metric{\D_{\Sigma_{m-1}}\log\rho_{m-1}}{\nu_m}
			+\sum_{k=2}^{m-1}H_{\Sigma_k}^2-\sum_{k=1}^{m-1}\lambda_k \label{eq-dim6:aux3}\\
			&\qquad -\sum_{k=1}^{m-1}\big(
				|h_{\Sigma_k}|^2+\Ric_{\Sigma_{k-1}}(\nu_k,\nu_k)
				+\metric{\D_{\Sigma_k}\log\rho_k}{\D_{\Sigma_k}\log v_k}
			\big) \nonumber
		\Big],
	\end{align}
	where $\lambda_k\geq0$ is the first eigenvalue of the stability operator for $\Sigma_k$. Define
	\[\Lambda=\sum_{k=1}^{m-1}\lambda_k,\quad \mathcal{G}=\sum_{k=1}^{m-1}\metric{\D_{\Sigma_k}\log\rho_k}{\D_{\Sigma_k}v_k},\quad \mathcal{R}=\sum_{k=1}^m\Ric_{\Sigma_{k-1}}(\nu_k,\nu_k)\]
	and
	\[\tilde{\mathcal{E}}=\sum_{k=1}^{m-1}|h_{\Sigma_k}|^2-\sum_{k=2}^{m-1}H_{\Sigma_k}^2.\]
	Note that $\Lambda,\mathcal{G},\mathcal{R}$ are defined in the same manner as in \cite[Lemma 3.4]{Brendle-Hirsch-Johne_2022}, while $\tilde{\mathcal{E}}$ differs from the expression $\mathcal{E}$ in \cite{Brendle-Hirsch-Johne_2022} by removing the bottom terms $k=m$. Rearranging the right hand side of \eqref{eq-dim6:aux3}, we obtain
	\begin{align}
		\frac{dH_\rho(t)}{dt}\int_{\Sigma_{m,t}}\varphi_t^{-1}
		&\leq \int_{\Sigma_{m,t}}\Big[
			-(\Lambda+\mathcal{G}+\tilde{\mathcal{E}}+\mathcal{R})-|h_{\Sigma_m}|^2
			-\metric{\D_{\Sigma_m}\log\varphi_t}{\D_{\Sigma_m}\log\tilde\rho_m} \nonumber\\
			&\qquad -H_{\Sigma_{m,t}}\metric{\D_{\Sigma_{m-1}}\log\rho_{m-1}}{\nu_m}
		\Big]. \label{eq-dim6:aux4}
	\end{align}
	In the estimate of gradient term $\mathcal{G}$ in \cite[Lemma 3.7]{Brendle-Hirsch-Johne_2022}, one used the condition that each $\Sigma_k$ has zero generalized mean curvature. However, through a careful look at the proof, we find that this fact is only used to convert the expression $\metric{\D_{\Sigma_k}\rho_k}{\nu_{k+1}}^2$ to $H_{\Sigma_{k+1}}^2$. In light of this, the proof of \cite[Lemma 3.7]{Brendle-Hirsch-Johne_2022} applies to our case to give
	\[\mathcal{G}\geq
		\sum_{k=2}^{m-1}\big(\frac12+\frac1{2(k-1)}\big)H_{\Sigma_k}^2
		+\frac m{2(m-1)}\Big(
			|\D_{\Sigma_m}\log\rho_{m-1}|^2
			+\metric{\D_{\Sigma_{m-1}}\log\rho_{m-1}}{\nu_m}^2
		\Big).\]
	The second term here appears in the end of the proof of \cite[Lemma 3.7]{Brendle-Hirsch-Johne_2022} but is discarded in the final statement. We need this extra positivity in our argument.
	
	The identity for $\mathcal{R}$ in \cite[Lemma 3.8]{Brendle-Hirsch-Johne_2022}, the grouping identities above \cite[Lemma 3.10]{Brendle-Hirsch-Johne_2022} for $1\leq k\leq m-1$, as well as the curvature inequalities for top and intermediate slices \cite[Lemma 3.11, 3.12]{Brendle-Hirsch-Johne_2022} remain applicable since they do not involve the second fundemantal form of the bottom slice $\Sigma_m$. As a result, we obtain the following inequality:
	\begin{equation}\label{eq-dim6:aux5}
		\mathcal{R}+\tilde{\mathcal{E}}+\mathcal{G}
		\geq \frac m{2(m-1)}\Big(
		|\D_{\Sigma_m}\log\rho_{m-1}|^2
		+\metric{\D_{\Sigma_{m-1}}\log\rho_{m-1}}{\nu_m}^2
		\Big).
	\end{equation}
	Combining \eqref{eq-dim6:aux4} and \eqref{eq-dim6:aux5} and using the estimate $|h_{\Sigma_m}|^2\geq H_{\Sigma_m}^2/(n-m)$, we have
	\[\begin{aligned}
		& \frac{dH_\rho(t)}{dt}\int_{\Sigma_{m,t}}\varphi_t^{-1} \\
		&\leq \int_{\Sigma_{m,t}}\Big[
		-|h_{\Sigma_m}|^2
		-\metric{\D_{\Sigma_m}\log\varphi_t}{\D_{\Sigma_m}\log\tilde\rho_m} -H_{\Sigma_m}\metric{\D_{\Sigma_{m-1}}\log\rho_{m-1}}{\nu_m} \\
		&\qquad\quad
		-\frac m{2(m-1)}\Big(
			|\D_{\Sigma_m}\log\rho_{m-1}|^2
			+\metric{\D_{\Sigma_{m-1}}\log\rho_{m-1}}{\nu_m}^2 \Big)
		\Big] \\
		&\leq \int_{\Sigma_{m,t}}\Big[ -\frac{H_{\Sigma_m}^2}{n-m} -H_{\Sigma_m}\metric{\D_{\Sigma_{m-1}}\log\rho_{m-1}}{\nu_m}
		-\frac m{2(m-1)}\metric{\D_{\Sigma_{m-1}}\log\rho_{m-1}}{\nu_m}^2 \\
		&\qquad\quad
		-\frac{m}{2(m-1)}|\D_{\Sigma_m}\log\rho_{m-1}|^2
		-\metric{\D_{\Sigma_m}\log\varphi_t}{\D_{\Sigma_m}(\log\rho_{m-1}+\log\varphi_t)}\Big].
	\end{aligned}\]
	The desired result $dH_\rho(t)/dt\leq0$ is obtained provided that the second to last and the last lines are both non-positive. Applying Young's inequality to the second to last line, we find that it is non-positive provided
	\begin{equation}\label{eq-dim6:aux2}
		\frac{2m}{(n-m)(m-1)}\geq1\ \Leftrightarrow\ n(m-1)\leq m^2+m,
	\end{equation}
	which is indeed true for $n=6$, $m\leq n-1$. The last line can be verified to be non-positive by Young's inequality again. As a result we have $dH_\rho(t)/dt\leq0$, which implies $ H_\rho(t)\leq0$ for $t>0$ combined with $H_\rho(0)=0$. This leads to a proof of the rigidity result, as explained at the beginning of this section.
\end{proof}

\vspace{12pt}

\noindent\textit{{Kai Xu,}}

\vspace{2pt}

\noindent\textit{{Department of Mathematics, Duke University, Durham, NC 27708, USA,}}

\vspace{2pt}

\noindent\textit{Email address: }\href{mailto:kx35@math.duke.edu}{kx35@math.duke.edu}.

\end{document}